\documentclass[a4paper,reqno]{amsart} 
\usepackage[utf8]{inputenc}                             
\usepackage{lmodern} 
\usepackage{microtype}
\usepackage{enumitem}
\usepackage{mathabx}

\newtheorem{theorem}{Theorem}[section]

\newtheorem{proposition}[theorem]{Proposition}
\newtheorem{lemma}[theorem]{Lemma}

\theoremstyle{definition}
\newtheorem{definition}[theorem]{Definition}
\newtheorem{example}[theorem]{Example}

\theoremstyle{remark}
\newtheorem{remark}[theorem]{Remark}
\numberwithin{equation}{section}

\begin{document}

\title{Convergence of ADAM for Lipschitz Objective Functions}

\author{Juan Ferrera}
\address{IMI, Departamento de An{\'a}lisis Matem{\'a}tico y
  Matemáticas Aplicadas,
  Facultad Ciencias Ma\-te\-m{\'a}ticas, Universidad Complutense, 28040, Madrid, Spain}
\email{ferrera@ucm.es}

\author{Javier G\'omez Gil}
\address{Departamento de An{\'a}lisis Matem{\'a}tico y
  Matemáticas Aplicadas,
  Facultad Ciencias Matem{\'a}ticas, Universidad Complutense, 28040, Madrid, Spain}
\email{gomezgil@ucm.es}
\thanks{The authors were partially supported by grant
 PGC2018-097286-B-100, Spanish  Agencia Estatal de Investigación,   http://dx.doi.org/10.13039/501100011033.}

\begin{abstract}
  The aim of this paper is to prove the exponential convergence, local
  and global, of Adam algorithm under precise conditions on the
  parameters, when the objective function lacks
  differentiability. More precisely, we require Lipschitz continuity,
  and control on the gradient whenever it exists. We provide also
  examples of interesting functions that satisfies the required
  restrictions.
\end{abstract}


\maketitle

\section{Introduction}

Adam (Adaptive moment estimation) is  one of the most used  optimization methods for training neural networks. 
It was introduced in \cite{KB}. In that paper
 the authors presented a result
on the algorithm convergence. 
However, there is a gap in the proof, and as a matter of fact there are counterexamples that prove that the result is false.
(See \cite{MR} and \cite{RKK}).
In \cite{BW}  a proof of the algorithm local convergence assuming strong differentiability restrictions on the cost function is given.
More recently, see \cite{DBB}, in a wider frame that includes also the zero momentum case, the authors present results for  non convex 
objective functions, but 
smoothness is also required. 

The aim of this paper is to remove the condition of the objective function differentiability. This is interesting since although it is usual to 
use functions like sum of squares ($\Vert\ \Vert_2^2$), non smooth objective functions as $\Vert\ \Vert_{\infty}^2$, could be interesting.

In this paper we will restrict to the deterministic case, and will provide local and global results that guarantee the quick convergence
assuming that the first order momentum parameter is small. First we introduce some notation and definitions.
Let $C:\mathbb{R}^N\to \mathbb{R}$ be a locally  Lipschitz function. It is well known
that for every $w\in \mathbb{R}^N$, Clarke's Generalized Gradient, $\partial C(w)\subset \mathbb{R}^N$, 
is a non empty convex closed set (see \cite{C}). It is also known, see \cite{F}, that
\begin{equation*}
  \partial_CC(w)=\text{co}\{ \lim_k\nabla C(w_k):\; f \, \text{is differentiable at }\, w_k, \; \lim_kw_k=w\}
\end{equation*}
(Remember that $C$ is differentiable a.e. by Rademacher Theorem).

Given $w\in \mathbb{R}^N$,  we define
$\zeta_w\in \partial_CC(w)$ as the nearest point at the origin of $\partial_CC(w)$. Observe that  convexity guarantee unicity of $\zeta_w$.
Nevertheless, if due to algorithm implementation requirements is better, there is no problem if we take
 $\zeta_w=\nabla C(w)$ whenever $C$ is differentiable at $w$.

We proceed to define the Generalized Adam algorithm:

\begin{definition} (Generalized Adam)
  Let  $\varepsilon >0$, $\beta_{1},\beta_{2}\in [0,1)$,
  $\alpha_n >0$ for every $n$.
  Let  $w_0, m_0, v_0 \in \mathbb{R}^N$ where the $v_0$ cordinates 
  are bigger or equal than $0$. 
  We define recursively 
  \begin{equation*}
    \begin{split}
      m_{n+1}=&\beta_{1} m_n+(1-\beta_{1})\zeta_{w_n},
      \\
      v_{n+1}=&\beta_{2} v_n+(1-\beta_{2})\zeta_{w_n}^2,
      \\ 
      w_{n+1}=&w_n-\alpha_n
    \frac{m_{n+1}}{\sqrt{v_{n+1}+(\varepsilon ,\dots ,\varepsilon )}}.
    \end{split}
  \end{equation*}
  where all the vector operations: we mean square, square root, product and division
  are defined component-wise.
\end{definition}

When  for a given $\alpha >0$,
$\alpha_n=\alpha \frac{\sqrt{1-\beta_2^{n+1}}}{1-\beta_1^{n+1}}$
for every $n$, we obtain the so called Adam algorithm. In this case Definition reads as follows

\begin{definition} (Adam)
  Let  $\varepsilon >0$, $\beta_{1},\beta_{2}\in [0,1)$,
  $\alpha >0$.
  Let  $w_0, m_0, v_0 \in \mathbb{R}^N$ where the $v_0$ coordinates 
  are bigger or equal than $0$. 
  We define recursively 
  \begin{equation*}
    \begin{split}
      m_{n+1}=&\beta_{1} m_n+(1-\beta_{1})\zeta_{w_n},
      \\
      v_{n+1}=&\beta_{2} v_n+(1-\beta_{2})\zeta_{w_n}^2,
      \\ 
      w_{n+1}=&w_n-\alpha \frac{\sqrt{1-\beta_2^{n+1}}}{1-\beta_1^{n+1}}
    \frac{m_{n+1}}{\sqrt{v_{n+1}+(\varepsilon ,\dots ,\varepsilon )}}.
    \end{split}
  \end{equation*}
 \end{definition}
 Of course,  if $f$ is $C^1$  both definitions agree with the usual one.
The following remark is immediate

\begin{remark}\label{formulasbasicas}
  \begin{equation*}
    m_{n+1}=\beta_1 ^{n+1}m_0+(1-\beta_1 )\sum_{i=0}^n \beta_1 ^{n-i}\zeta_{w_i},
  \end{equation*}
  \begin{equation*}
    v_{n+1}=\beta_2 ^{n+1}v_0+(1-\beta_2 )\sum_{i=0}^n \beta_2 ^{n-i}\zeta_{w_i}^2.
  \end{equation*}
\end{remark}

By the moment we will focus on the particular case. We may describe the algorithm
as a non autonomous discrete dynamical system as follows:
If $X= \mathbb{R}^N\times [0,+\infty)^{N}\times \mathbb{R}^{N}$, 
for $x=(m,v,w)\in X$ we define
\begin{equation*}
  \begin{split}
    M(m,w)=&\beta_1m+(1-\beta_1)\zeta_w,
    \\
    V(v,w)=&\beta_2v+(1-\beta_2)\zeta_w^2,
  \end{split}
\end{equation*}
$\Gamma :X\to X$ as
\begin{equation*}
  \Gamma (x)=
  \begin{pmatrix}
    M(m,w) \\ V(v,w) \\[0.1cm]
    w-\alpha
    \frac{M(m,w)}{\sqrt{V(v,w)+(\varepsilon ,\dots ,\varepsilon )}}
  \end{pmatrix}
\end{equation*}
$\Omega :\mathbb{N}\times X\to X$ as
\begin{equation*}
  \Omega (n,x)=
  \begin{pmatrix}
    0 \\0 \\[0.1cm]
    -\alpha \biggl( \frac{\sqrt{1-\beta_2^{n+1}}}{1-\beta_1^{n+1}}-1\biggr)
    \frac{M(m,w)}{\sqrt{V(v,w)+(\varepsilon ,\dots ,\varepsilon )}}
  \end{pmatrix}
\end{equation*}
and finaly 
\begin{equation*}
  \Theta:\mathbb{N}\times X\longrightarrow X
\end{equation*}
as $\Theta(n,x)=\Gamma (x)+\Omega (n,x)$. 
In this way we have: $x_{n+1}=\Theta(n,x_n)$. This notation allow us to establish
the following result:

\begin{theorem}\label{sec:plant-del-probl}
 If $\zeta_{w^*}=0$ then $x^*=(0,0,w^*)$
  is a fixed point of $\Gamma$, and 
  $\Omega (n,x^*)=0$ for every $n\in \mathbb{N}$,
   and consequently  $\Theta(n,x^*)=x^*$ for every $n\in \mathbb{N}$.
\end{theorem}
In this paper, unless we say otherwise,
 $\Vert \ \Vert$ will denote the euclidean norm on  $\mathbb{R}^N$.

\section{Local Convergence}

The aim of this section is to study the algorithm asymptotic convergence
 near to critical points of the cost function $C$. 
 In order to proceed, we will generalize the following version of Corollary III.4 of \cite{BW}:

 \begin{theorem}\label{BWtheorem}
 Assume that $C$ is a $C^2$ function that attains a minimum at $w^*$.  Let us denote the eigenvalues of the  
 (positively defined)
 Hessian matrix of $C$ at $w^*$ by $\lambda_i$, $i=1,\dots ,N$. Assume that
 \begin{equation*}
   \frac{\alpha (1-\beta_1)}{\sqrt{\varepsilon}}\max_i\{ \lambda_i\}<2\beta_1+2.
 \end{equation*}
 Then there exists a neighborhood of $x^*=(0,0,w^*)\in \mathbb{R}^{3N}$, $V$, such that 
 if $x_0=(0,0,w_0)\in V$ then $(x_n)_n$ converges exponentially to $x^*$.
  \end{theorem}
  
  Specifically, we will remove the fact that $C$ is $C^2$,
  although  we must 
 modify the parameters restriction, and we will require a control on $\nabla C(w)$ whenever it exists.
 Indeed,  we will impose only conditions on 
 $\nabla C(w)$ where $C$ is differentiable, and first of all we see how $\zeta_w$
 inherits these properties.

 \begin{lemma}\label{L1}
  Let $U$ be a neighborhood of $w^{*}$.
  \begin{enumerate}[label=\textup{(\arabic*)},wide]
  \item \label{item:1} If $\Vert \nabla C(w)\Vert \leq \mu \Vert
    w-w^*\Vert $ for every $w\in U$ where $C$ is
    differentiable then
    $\Vert \zeta_w\Vert \leq \mu \Vert w-w^*\Vert $ for every $w\in U$, and consequently
    $\zeta_{w^*}=0$.
  \item \label{item:7} If $\langle \nabla C(w),w-w^*\rangle \geq \delta \Vert w-w^*\Vert^2$
       for every $w$ where $C$ is differentiable,
       then $\langle \zeta_w,w-w^*\rangle \geq \delta \Vert w-w^*\Vert^2$ for every $w\in U$.
  
  \end{enumerate}
\end{lemma}
\begin{proof}
    If $\zeta_w=\lim_k\nabla C(w_k)$, taking limits the result follows immediately. In the general case
    $\zeta_w=\lambda_1\zeta_1+\dots +\lambda_m\zeta_m$ with $\lambda_i>0$ and $\sum_i\lambda_i=1$
    and $\zeta_i$ satisfying the corresponding inequality. Hence
    \begin{enumerate}
  \item $\Vert \zeta_w\Vert \leq \sum_i \lambda_i \Vert \zeta_i \Vert 
  \leq \mu \Vert w-w^*\Vert $ 
\end{enumerate}
and
\begin{enumerate}[resume]
  \item $\langle \zeta_w,w-w^*\rangle =\sum_i\lambda_i\langle \zeta_i,w-w^*\rangle 
  \geq \delta \Vert w-w^*\Vert^2$.
\end{enumerate}
\end{proof}

We will take advantage of the second condition in a slightly different and apparently stronger
form.

\begin{lemma}\label{sec:convergencia-local-2}
Let $U$ be a neighborhood  of $w^{*}$. We denote $\mathbf{a}=(a,\dotsc,a)$, $a>0$. 
If for $0<\delta \leq \mu$ we have  that: 
  \begin{enumerate}[label=\textup{(\arabic*)}]
  \item \label{item:28} $\Vert \zeta_w\Vert \leq \mu \Vert w-w^*\Vert
    $
  \end{enumerate}
  and
  \begin{enumerate}[resume*]
  \item\label{item:29}
    $\displaystyle\langle \zeta_w,w-w^*\rangle \geq
    \delta \Vert w-w^*\Vert^2$
  \end{enumerate}
   for every $w\in U$, then 
    \begin{equation}
    \langle \mathbf{b}\zeta_w,w-w^*\rangle \geq (\max_{i} b_{i})
   \frac{3 \delta}{4} \Vert w-w^*\Vert^2\label{eq:9}
  \end{equation}
 for every $\mathbf{b}\in \mathbb{R}^{N}$ satisfying $\Vert \mathbf{b}-\mathbf{a}\Vert_\infty
  <\eta=\dfrac{a\delta}{4\mu+3\delta}$, and every $w\in U$.
\end{lemma}
\begin{proof}
  By $(1)$
  \begin{equation*}
    \langle (\mathbf{b}-\mathbf{a})\zeta_w,w-w^*\rangle=\langle \zeta_w,(\mathbf{b}-\mathbf{a})(w-w^*)\rangle 
     \leq \mu \Vert w-w^*\Vert^{2}\,\Vert \mathbf{b}-\mathbf{a}\Vert_{\infty}
  \end{equation*}
  hence, by $(2)$,
  \begin{align*}
    \langle \mathbf{b}\zeta_w,w-w^*\rangle= & \langle (\mathbf{b}-\mathbf{a})\zeta_w,w-w^*\rangle+
    \langle \mathbf{a}\zeta_w,w-w^*\rangle
    \\ \geq & a\delta \Vert w-w^*\Vert^2- \mu
              \Vert w-w^*\Vert^{2}\,\Vert \mathbf{b}-\mathbf{a}\Vert_{\infty} \\= &
                                                                           \left(a\delta-\mu \Vert \mathbf{b}-\mathbf{a}\Vert_{\infty} \right)\Vert w-w^*\Vert^2
                                                                           > \left(a\delta-\mu \eta \right)\Vert w-w^*\Vert^2.
  \end{align*}
  As 
  \begin{equation*}
    a\delta-\mu \eta = a\delta\left(1-
      \frac{\mu}{4\mu+3\delta} \right)=3a\delta \frac{\mu+\delta}{4\mu+3\delta}
      =(a+\eta)\frac{3\delta}{4}
  \end{equation*}
 since 
  \begin{equation*}
  a \frac{\mu+\delta}{4\mu+3\delta}=a\frac14\left[ 1+\frac{\delta}{4\mu+3\delta}\right] =
  \frac14 (a+\eta)
  \end{equation*}
  the fact that $a+\eta \geq (\max_{i} b_{i})$ give us the result.
\end{proof}

In order to study the behavior of $\Theta$, we begin with the autonomous part
$\Gamma$. The components $\Gamma_1$ and $\Gamma_2$ behave simply.

\begin{proposition}\label{sec:convergencia-local}
  Assume that $C:\mathbb{R}^N\to \mathbb{R}$ is locally  Lipschitz  and
  there exist $\mu >0$ and   $R\in (0,1)$ such that 
  $\Vert \nabla C(w)\Vert \leq \mu \Vert w-w^*\Vert $ provided that
  $C$ is differentiable at $w$ and
  $w\in B(w^*,R)$.
  Then
  \begin{enumerate}[label=\textup{(\arabic*)}]
  \item \label{item:5} $\displaystyle
    \Vert \Gamma_1(x)\Vert \leq
    \beta_1\Vert m\Vert +\mu (1-\beta_1)\Vert w-w^*\Vert
    $
  \item \label{item:6} $\displaystyle
    \Vert \Gamma_2(x)\Vert \leq
    \beta_2\Vert v\Vert +\mu ^2(1-\beta_2)\Vert w-w^*\Vert
    $
  \end{enumerate}
  for every $x=(m,v,w)\in  X$ satisfying $\Vert w-w^*\Vert <R$.
\end{proposition}
\begin{proof}
For the second part it is enough to observe that 
$\Vert \zeta_w^2\Vert \leq \Vert \zeta_w\Vert^2$.
\end{proof}

The third component requires some more work.

\begin{proposition}
  \label{P1}
  Assume that $C:\mathbb{R}^N\rightarrow \mathbb{R}$ is locally Lipschitz, and 
  there exist $0<\delta \leq \mu$ and $R\in (0,1)$ such that:
  \begin{enumerate}[label=\textup{(\arabic*)}]
  \item \label{item:30} $\Vert \nabla C(w)\Vert \leq \mu \Vert w-w^*\Vert $
  \item \label{item:31} $\displaystyle\langle \nabla C(w),w-w^*\rangle \geq
    \delta \Vert w-w^*\Vert^2$
  \end{enumerate}
  provided that $C$ is differentiable at $w\in B(w^*,R)$.
  Then, for every $w\in B(w^*,R)$,  if 
  $\Vert V(v,w)\Vert_{\infty}<\left[(1+\frac{\delta}{2\mu+\delta})^2-1\right]\varepsilon$,
we have that  
  \begin{equation*}
    \Vert \Gamma_3(m,v,w)-w^*\Vert \leq \frac{\alpha \beta_1}{\sqrt{\varepsilon}}\Vert m\Vert +D \Vert w-w^*\Vert
  \end{equation*}
  where
  \begin{equation*}
    D=
    \Biggl( 1- \frac{9\alpha \delta(1-\beta_1)}{8\sqrt\varepsilon}  +
    \frac{\alpha^2 (1-\beta_1)^2\mu^2}{\varepsilon}\Biggr) ^{\frac{1}{2}}.
  \end{equation*}
\end{proposition}
\begin{proof}
  Remember that we denote $(m,v,w)=x$ and that $v\geq 0$. We have that
  \begin{equation*}
    \Vert \Gamma_3(x)-w^*\Vert  \leq
    \frac{\alpha \beta_1}{\sqrt{\varepsilon}}\Vert m\Vert +
    \biggl\Vert
    (w-w^*)-\frac{\alpha (1-\beta_1)\zeta_w}{\sqrt{V(v,w)+(\varepsilon ,\dots ,\varepsilon )}}
    \biggr\Vert.
  \end{equation*}
  We estimate the second term. If we denote
  \begin{equation*}
    \mathbf{b}=(b_1,\dotsc,b_N)=
    \frac{\alpha (1-\beta_1)}{\sqrt{V(v,w)+(\varepsilon ,\dots ,\varepsilon )}},
  \end{equation*}
  it is immediate that
  \begin{equation*}
   \frac{\alpha (1-\beta_1)}{(1+\frac{\delta}{2\mu+\delta})\sqrt{\varepsilon}}\leq
    b_{i}\leq \frac{\alpha (1-\beta_1)}{\sqrt\varepsilon}
  \end{equation*}
  therefore, if we denote
  \begin{equation*}
    a=\frac{\alpha (1-\beta_1)}{2}\left[\frac1{\sqrt\varepsilon}  +
     \frac{1}{(1+\frac{\delta}{2\mu+\delta})\sqrt{\varepsilon}}\right]= \frac{\alpha (1-\beta_1)}{2\sqrt\varepsilon}\, \frac{4\mu+3\delta}{2\mu+2\delta}
  \end{equation*}
  we have that
  \begin{align*}
    |b_{i}-a|\leq & \frac{\alpha (1-\beta_1)}{2}\left[\frac1{\sqrt\varepsilon}  -
       \frac{1}{(1+\frac{\delta}{2\mu+\delta})\sqrt{\varepsilon}}\right]
                    =\frac{\alpha
                    (1-\beta_1)}{2\sqrt\varepsilon}\,\frac{\delta}{2\mu+2\delta}
    \\ =&
          \frac{\alpha (1-\beta_1)}{2\sqrt\varepsilon}\,\frac{4\mu+3\delta}{2\mu+2\delta}\,
          \frac{\delta}{4\mu+3\delta}=  \frac{a\delta}{4\mu+3\delta}
  \end{align*}
  for every  $i=1,\dotsc,N$. Hence, invoking Lemma
  \ref{sec:convergencia-local-2} we get
  \begin{equation*}
    \langle \mathbf{b}\zeta_w,w-w^*\rangle \geq (\max_{i} b_{i})\frac{3\delta}{4} \Vert w-w^*\Vert^2
  \end{equation*}
  and consequently
  \begin{align*}
      \Vert (w-w^*)-&\mathbf{b}\zeta_w\Vert^2= \Vert w-w^*\Vert^2-2\langle \mathbf{b}\zeta_w,w-w^*\rangle +
                      \Vert \mathbf{b}\zeta_w\Vert^2 \\ \leq &
                                                       \Vert w-w^*\Vert^2-2\langle \mathbf{b}\zeta_w,w-w^*\rangle +
                                                       \frac{\alpha^2 (1-\beta_1)^2}{\varepsilon}\Vert \zeta_w\Vert^2 \\
      \leq &
             \Vert w-w^*\Vert^2-(\max_i b_{i})\frac{3\delta}{2}  \Vert w-w^*\Vert^2+
             \frac{\alpha^2 (1-\beta_1)^2}{\varepsilon}\Vert \zeta_w\Vert^2 \\
      \leq &
             \Vert w-w^*\Vert^2 \left(1 -(\max_i b_{i})\frac{3\delta}{2}  +
             \frac{\alpha^2 (1-\beta_1)^2\mu^2}{\varepsilon}\right) 
  \end{align*}
  by Lema \ref{L1}. As 
    \begin{equation*}
\max_i b_{i}\geq  \frac{2\mu+\delta}{2\mu+2\delta}\,\frac{\alpha
  (1-\beta_1)}{\sqrt{\varepsilon}}\geq  \,\frac{3\alpha
  (1-\beta_1)}{4\sqrt{\varepsilon}}
\end{equation*}
we deduce that
  \begin{equation*}
    \Vert \Gamma_3(x)-w^*\Vert \leq
    \frac{\alpha \beta_1}{\sqrt{\varepsilon}}\Vert m\Vert +
    D\Vert w-w^*\Vert.
  \end{equation*}
\end{proof}

In order to join Propositions \ref{sec:convergencia-local} and
\ref{P1}
it is convenient to consider in $\mathbb{R}^{3N}$ the norm
\begin{equation*}
      \Vert x\Vert _{\infty} =\max\{\Vert m\Vert ,\Vert v\Vert ,\Vert w\Vert\}.
    \end{equation*}
Now we may describe the behavior of $\Gamma$

\begin{theorem}
  \label{COR}
  Assume that $C:\mathbb{R}^N\rightarrow \mathbb{R}$ is locally Lipschitz, and 
  there exist $0<\delta \leq \mu$ and $R\in (0,1)$ such that:
  \begin{enumerate}[label=\textup{(\arabic*)}]
  \item \label{item:301} $\Vert \nabla C(w)\Vert \leq \mu \Vert w-w^*\Vert $
  \item \label{item:311} $\displaystyle\langle \nabla C(w),w-w^*\rangle \geq
    \delta \Vert w-w^*\Vert^2$
  \end{enumerate}
  provided that $C$ is differentiable at $w\in B(w^*,R)$.

Then there exists  $r\in (0,R)$ such that if we denote  $x^{*}=(0,0,w^{*})$,
  \begin{equation*}
    \Vert \Gamma (x)-x^*\Vert _{\infty}\leq K\Vert x-x^*\Vert _{\infty}
  \end{equation*}
  for every $x\in B_{\infty}(x^*,r)\cap X$, 
  where
  \begin{equation*}
    K=\max \bigg(
    \beta_1+\mu (1-\beta_1),
    \beta_2+\mu ^2(1-\beta_2), \frac{\alpha\beta_{1}}{\sqrt{\varepsilon}}+D\bigg)
  \end{equation*}
  and
  \begin{equation}\label{eq:11}
    D=\bigg( 1- \frac{9\alpha \delta (1-\beta_1)}{8\sqrt{\varepsilon}} +
    \frac{\alpha^2 (1-\beta_1)^2\mu^2}{\varepsilon}\bigg) ^{\frac{1}{2}}.
  \end{equation}
\end{theorem}
\begin{proof}
 First, we observe that
  \begin{equation*}
    \begin{split}
      \Vert \Gamma_1(x)\Vert &\leq
                               \Big( \beta_1+\mu (1-\beta_1)\Big) \Vert x-x^*\Vert_{\infty},
      \\
      \Vert \Gamma_2(x)\Vert &\leq
                               \Big( \beta_2+\mu ^2(1-\beta_2)\Big) \Vert x-x^*\Vert_\infty.
    \end{split}
  \end{equation*}
  In order to estimate the third component
  it is enough to take
  \begin{equation*}
  r=\frac{\varepsilon}{(\beta_2+(1-\beta_2)\mu ^2)}\left[\left(1+\frac{\delta}{2\mu+\delta}\right)^2-1\right]
\end{equation*}
since from 
\begin{equation*}
  \Vert V\Vert _{\infty}\leq  
  \Big( \beta_2+(1-\beta_2)\mu ^2\Big) \Vert x-x^*\Vert _{\infty}
\end{equation*}
it follows that  $\Vert V\Vert _{\infty}\leq
  \left[(1+\frac{\delta}{2\mu+\delta})^2-1\right]\varepsilon$,
  provided that $\Vert x-x^*\Vert _{\infty}<r$. This inequality allows us to invoke
  Proposition \ref{P1}, and deduce
  \begin{equation*}
    \Vert \Gamma_3(x)-w^*\Vert \leq
    \left(\frac{\alpha\beta_{1}}{\sqrt{\varepsilon}}+D\right) \Vert x-x^*\Vert _{\infty}.
  \end{equation*}
\end{proof}

It is obvious that the radicand that defines $D$ is positive, so $D$ is a well defined positive number,
but moreover it will be useful that $D<1$; this is equivalent to
\begin{equation}
1-\beta_1<\frac{9\delta \sqrt{\varepsilon}}{8\mu^2\alpha}.
\end{equation}
Anyway, next Lemma gives us an effective bound,
under a restriction on the parameters,  which is enough for our goals.

\begin{lemma}\label{L3}
  If
  \begin{equation}
    1-\beta_1=\frac{\delta
      \sqrt\varepsilon}{2\alpha \mu^2}\label{eq:8}
  \end{equation}
  then
  \begin{equation*}
   D = \bigg( 1-\frac98\,\frac{\alpha (1-\beta_1)}{\sqrt{\varepsilon}}\delta  +
    \frac{\alpha^2 (1-\beta_1)^2\mu^2}{\varepsilon}\bigg) ^{\frac{1}{2}}<1-\frac{\delta^2}{8\mu^2}.
  \end{equation*}
\end{lemma}

\begin{proof}
  \begin{equation*}
      D = \bigg( 1-\frac9{16}\frac{\delta^2}{\mu^2}+\frac14                                      \frac{\delta^2}{\mu^2} \bigg) ^{\frac{1}{2}} =
\bigg( 1-\frac5{16}\frac{\delta^2}{\mu^2}\bigg) ^{\frac{1}{2}}<
1-\frac5{32}\frac{\delta^2}{\mu^2}< 1-\frac{\delta^2}{8\mu^2}
  \end{equation*}
since $\sqrt{1-y}<1-\frac{y}{2}$  provided that $0<y<1$.
\end{proof}

In order to prove that $\Gamma$ is contractive at $x^*$, we introduce a new norm on $\mathbb{R}^{3N}$. For $A\geq 1$
we define
 \begin{equation*}
      \vvvert x\vvvert =\max\{\Vert m\Vert ,\Vert v\Vert ,A\Vert w\Vert\}.
    \end{equation*}
 
  \begin{theorem}\label{T1}
    Assume that $C:\mathbb{R}^N\rightarrow \mathbb{R}$ is locally Lipschitz, and 
  there exist $0<\delta \leq \mu$ and $R\in (0,1)$ such that:
  \begin{enumerate}[label=\textup{(\arabic*)}]
  \item \label{item:302} $\Vert \nabla C(w)\Vert \leq \mu \Vert w-w^*\Vert $
  \item \label{item:312} $\displaystyle\langle \nabla C(w),w-w^*\rangle \geq
    \delta \Vert w-w^*\Vert^2$
  \end{enumerate}
  provided that $C$ is differentiable at $w\in B(w^*,R)$. Let $A>\mu$ (if $\mu <1$ we take $A=1$).
If
 \begin{enumerate}[label=\textup{(\alph*)}]
    \item \label{item:25}
      $\displaystyle
     \frac{\delta\sqrt\varepsilon}{2\mu^{2}}<
      \alpha\leq
      \frac{\delta\sqrt\varepsilon}{2\mu^{2}}\left(1+\frac{\delta}{4A}\right)
      $ 
    \end{enumerate}
    and
    \begin{enumerate}[resume*]
    \item \label{item:26}
      $\displaystyle
     1-\beta_{1}= \frac{\delta \sqrt\varepsilon}{2\alpha\mu^2}$,
    \end{enumerate}
    then there exist $r\in (0,R)$ and  $L_0<1$ such that
    $\vvvert \Gamma (x)-x^*\vvvert \leq L_0\vvvert x-x^*\vvvert$
    for every $x\in B_{\infty}(x^*,r)\cap X$.
  \end{theorem}

\begin{proof}
 We take $r$ like in Theorem \ref{COR} proof's.
    If $\Vert x-x^{*}\Vert<r$ then
    \begin{align*}
      \vvvert\Gamma (x)-x^*\vvvert  \leq &
                                           \max\biggl\{\beta_1\Vert m\Vert +\mu (1-\beta_1)\Vert w-w^*\Vert  ,
                                           \beta_2\Vert v\Vert  \\ & +\mu ^2(1-\beta_2)\Vert w-w^*\Vert+
                                                                     A\Bigl( \frac{\alpha \beta_1}{\sqrt{\varepsilon}}\Vert m\Vert +
                                                                     D\Vert w-w^*\Vert
                                                                     \Bigr) \biggr\}\\ \leq &
                                                                                              \max\biggl\{\beta_{1}+\frac{\mu}{A}(1-\beta_{1}),
                                                                                              \beta_{2}+\frac{\mu^{2}}{A^{2}}(1-\beta_{2}), A \frac{\alpha \beta_1}{\sqrt{\varepsilon}} +
                                                                                              D \biggr\}\vvvert x-x^*\vvvert.
    \end{align*}
    We have to check that
    \begin{equation*}
      L_{0}= \max\biggl\{\beta_{1}+\frac{\mu}{A}(1-\beta_{1}),
      \beta_{2}+\frac{\mu^{2}}{A^{2}}(1-\beta_{2}), A \frac{\alpha \beta_1}{\sqrt{\varepsilon}} +
      D \biggr\} <1.
    \end{equation*}
    As $\mu <A$, the first two terms in the max are smaller than $1$.
     Last, by  \ref{item:26} y \ref{item:25}
    \begin{equation*}
      A\frac{\alpha\beta_{1}}{\sqrt\varepsilon}=A\left[ \frac{\alpha}{\sqrt\varepsilon}-
        \frac{\alpha(1-\beta_{1})}{\sqrt\varepsilon}\right]= A
      \left[\frac{\alpha}{\sqrt\varepsilon}-\frac{\delta}{2\mu^{2}}\right]\leq A\frac{\delta}{2\mu^{2}}  \frac{\delta}{4A}=\frac{\delta^{2}}{8\mu^{2}}
    \end{equation*}
    hence $D+ A\frac{\alpha\beta_{1}}{\sqrt\varepsilon}<1$ by Lemma \ref{L3} .
 \end{proof}

As we observe in the Theorems's statement, if
 $\mu<1$ then $\vvvert\,\vvvert=\Vert \,\Vert_\infty$ works.

\begin{remark}
   Observe  that, as
\begin{equation*}
  \left(1+\frac{\delta}{4A}\right)  ^{-1}     \leq \frac{\delta  \sqrt\varepsilon}{2\alpha \mu^{2}}
\end{equation*}
we have
\begin{equation*}
  \beta_1\leq 1-\left(1+\frac{\delta}{4A}\right) ^{-1}\leq 1-\left(1+\frac{1}{4}\right) ^{-1}
  = \frac{1}{5}
\end{equation*}
since $\delta <A$.
\end{remark}

Now we proceed with the non autonomous term $\Omega$.
  
\begin{theorem}\label{T2}
 Assume that $C:\mathbb{R}^N\rightarrow \mathbb{R}$ is locally Lipschitz, and 
  there exist $0<\delta \leq \mu$ and $R\in (0,1)$ such that:
  \begin{enumerate}[label=\textup{(\arabic*)}]
  \item \label{item:303} $\Vert \nabla C(w)\Vert \leq \mu \Vert w-w^*\Vert $
  \item \label{item:313} $\displaystyle\langle \nabla C(w),w-w^*\rangle \geq
    \delta \Vert w-w^*\Vert^2$
  \end{enumerate}
  provided that $C$ is differentiable at $w\in B(w^*,R)$.
  Then for every $\beta_1,\beta_2\in (0,1)$,
   there exist $K_0>0$ and  $\beta \in (0,1)$ such that
  \begin{equation*}
    \Vert \Omega (n,x)\Vert _{\infty}\leq K_0\beta ^n\Vert x-x^*\Vert _{\infty}.
  \end{equation*}
\end{theorem}
\begin{proof} We denote  
$K_1=\frac{\alpha}{\sqrt{\varepsilon}} \big( \beta_1+\mu (1-\beta_1)\big)$, 
$K_2=(1-\beta_1)^{-1}\Big( (1-\beta_1)+\sqrt{1-\beta_2}\Big)^{-1}$,
$K_0=2K_1K_2$ and $\beta =\max \{ \beta_1,\beta_2\}<1$. Then

  \begin{equation*}
    \begin{split}
      \Vert \Omega (n,x)\Vert _{\infty}= &
                                           \left\Vert
                                           -\alpha \bigg( \frac{\sqrt{1-\beta_2^{n+1}}}{1-\beta_1^{n+1}}-1\bigg)
                                           \frac{M(m,w)}{\sqrt{V(v,w)+(\varepsilon ,\dots ,\varepsilon )}}
                                           \right\Vert
      \\ =&
            \left\Vert
            -\alpha \bigg( \frac{\sqrt{1-\beta_2^{n+1}}}{1-\beta_1^{n+1}}-1\bigg)
            \frac{\beta_1m+(1-\beta_1)\zeta_w}{\sqrt{V(v,w)+(\varepsilon ,\dots ,\varepsilon )}}
            \right\Vert
      \\ \leq&
               \frac{\alpha}{\sqrt{\varepsilon}}
               \bigg| \frac{\sqrt{1-\beta_2^{n+1}}}{1-\beta_1^{n+1}}-1\bigg|
               \Vert \beta_1m+(1-\beta_1)\zeta_w\Vert
      \\ \leq&
               \bigg| \frac{\sqrt{1-\beta_2^{n+1}}}{1-\beta_1^{n+1}}-1\bigg|
               \frac{\alpha}{\sqrt{\varepsilon}}
               \big( \beta_1+\mu (1-\beta_1)\big) \Vert x-x^*\Vert _{\infty}
      \\ =&
               \bigg| \frac{\sqrt{1-\beta_2^{n+1}}}{1-\beta_1^{n+1}}-1\bigg|
               K_1 \Vert x-x^*\Vert _{\infty}
      \\ =&
      \frac{1}{(1-\beta_1^{n+1})}
            \bigg| \frac{1-\beta_2^{n+1}-(1-\beta_1^{n+1})^2}{(1-\beta_1^{n+1})+\sqrt{1-\beta_2^{n+1}}}
            \bigg|  K_1 \Vert x-x^*\Vert _{\infty}
      \\  \leq &
                 \frac{1}{(1-\beta_1)}
                 \bigg|
                 \frac{1-\beta_2^{n+1}-(1-\beta_1^{n+1})^2}{(1-\beta_1)+\sqrt{1-\beta_2}}
                 \bigg| K_1 \Vert x-x^*\Vert _{\infty}
       \\ =&
                 K_1K_2\big| 1-\beta_2^{n+1}-(1-\beta_1^{n+1})^2\big| \Vert x-x^*\Vert _{\infty}
         \\=&
                 K_1 K_2\big| 2\beta_1^{n+1}-\beta_2^{n+1}-\beta_1^{2(n+1)}\big| \Vert x-x^*\Vert _{\infty}
                 \leq
                                                                                          2K_1 K_2\beta ^{n+1}\Vert x-x^*\Vert _{\infty}
 \end{split}
  \end{equation*}
since $\beta_1,\beta_2\in (0,1)$.
\end{proof}

Observe that for the norm $\vvvert\ \vvvert$ we have
\begin{equation}
\vvvert \Omega (n,x) \vvvert =
A\Vert \Omega (n,x)\Vert _{\infty}\leq AK_0\beta ^n\Vert x-x^*\Vert _{\infty} 
\leq AK_0\beta^n\vvvert x-x^*\vvvert.
\end{equation}

Finally, we may establish the main result of this section:

\begin{theorem}\label{mainlocal}
   Assume that $C:\mathbb{R}^N\rightarrow \mathbb{R}$ is locally Lipschitz, and 
  there exist $0<\delta \leq \mu$ and $R\in (0,1)$ such that:
  \begin{enumerate}[label=\textup{(\arabic*)}]
  \item \label{item:304} $\Vert \nabla C(w)\Vert \leq \mu \Vert w-w^*\Vert $
  \item \label{item:314} $\displaystyle\langle \nabla C(w),w-w^*\rangle \geq
    \delta \Vert w-w^*\Vert^2$
  \end{enumerate}
  provided that $w\in B(w^*,R)$ and $C$ is differentiable at $w$. Let $A>\mu$ (if $\mu <1$ we take $A=1$).
If
  \begin{enumerate}[label=\textup{(\alph*)}]
  \item \label{item:11}
      $\displaystyle 
     \frac{\delta\sqrt\varepsilon}{2\mu^{2}}<
      \alpha\leq
      \frac{\delta\sqrt\varepsilon}{2\mu^{2}}\left(1+\frac{\delta}{4A}\right)
      $  
    \end{enumerate}
    and
 \begin{enumerate}[resume*]     
  \item \label{item:12}
    $\displaystyle
    \beta_{1} =1-\frac{\delta \sqrt{\varepsilon}}{2\alpha\mu^2}$,
  \end{enumerate}
  then there exists a neighborhood 
  $U$ of $x^*$ such that $(x_n)_n$ converges exponentially to  $x^*$ for every 
 $x_0\in U$. In other words, there exist $n_0$ , $K>0$,  and $L\in (0,1)$
 such that if $x_0\in U$ then
  \begin{equation*}
    \vvvert x_n-x^*\vvvert \leq KL^{n-n_0}\vvvert x_0-x^*\vvvert  \quad \textrm{for every}\ n\geq n_0.
  \end{equation*}
\end{theorem}

\begin{proof}
By Theorems \ref{T1} y \ref{T2},
 there exist $U_0=B(x^*,r)\cap X$, $L_0<1$, $\beta <1$ and $K_0>0$  such that
 $\vvvert\Gamma (x)-x^*\vvvert\leq L_0\vvvert x-x^*\vvvert$
  and $\vvvert\Omega (n,x)\vvvert\leq K_0\beta ^n\vvvert x-x^*\vvvert $ for every  $x\in U_0$.
  
  Let us denote $L=\frac{L_0+1}{2}\in (L_0,1)$.
Let $n_0$ satisfies $L_0+K_0\beta ^{n_0}<L$. We define
  \begin{equation*}
    K=\max\biggl\{1,\prod_{i=0}^n(L_0+K_0\beta ^{i}):\; 0\leq n\leq n_0-1\biggr\}.
  \end{equation*}
  The neighborhood of $x^*$ will be $U=B(x^*,\frac{r}{K})\cap X$. We claim that
  $x_n\in U_0$ for every $n\geq 0$ provided that $x_0\in U$.
  
  We will prove
   \begin{equation*}
    \vvvert x_n-x^*\vvvert \leq KL^{n-n_0}\vvvert x_0-x^*\vvvert  \quad \textrm{for every}\ n\geq n_0
  \end{equation*}
by induction. For $n=n_0$ we have
 \begin{equation*}
    \vvvert x_{n_0}-x^*\vvvert  =\vvvert \Theta (n_0-1,x_{n_0-1})-x^*\vvvert  \leq
                                    \vvvert \Gamma (x_{n_0-1})-x^*\vvvert
                                    +\vvvert \Omega (n_0-1,x_{n_0-1})\vvvert 
 \end{equation*}
 \begin{equation*}
 \leq (L_0+K_0\beta ^{n_0-1})\vvvert x_{n_0-1}-x^*\vvvert \leq \dots \leq
 K\vvvert x_{0}-x^*\vvvert.
 \end{equation*}
Assume that the inequality holds for $x_{n_0},\dots , x_n$, then
 \begin{equation*}
    \begin{split}
      \vvvert x_{n+1}-x^*\vvvert \leq &(L_0+K_0\beta ^n)\vvvert x_n-x^*\vvvert
       \leq 
                (L_0+K_0\beta ^n)KL^{n-n_0}\vvvert x_0-x^*\vvvert \\ \leq &
                                                                            LKL^{n-n_0}\vvvert x_0-x^*\vvvert =KL^{n+1-t_0}\vvvert x_0-x^*\vvvert.
    \end{split}
  \end{equation*}
  
  In order to prove the claim we will proceed by induction again. For $n=0$ it is trivial,
 assume that $x_0,\dots ,x_n\in U$, as we see above
  \begin{equation*}
      \vvvert x_{n+1}-x^*\vvvert \leq (L_0+K_0\beta ^n)\vvvert x_n-x^*\vvvert.
  \end{equation*}
If $n\geq n_0$ then $(L_0+K_0\beta ^n)<L<1$, hence $\vvvert x_{n+1}-x^*\vvvert<\vvvert x_n-x^*\vvvert$
and we are done.  Otherwise, if $n<n_0$, repeating the argument we obtain
  \begin{equation*}
    \vvvert x_{n+1}-x^*\vvvert\leq \prod_{i=0}^n(L_0+K_0\beta ^{i})
    \vvvert x_0-x^*\vvvert
    \leq K\vvvert x_0-x^*\vvvert<r.
  \end{equation*}
  \end{proof}

%
%
%
%
%
%
Next remark allows us to observe how differentiability conditions on $C$ in Theorem \ref{BWtheorem}
are stronger than Theorem \ref{mainlocal}  ones: let us see it.

\begin{remark}
Assume that $C$ is $C^2$,  $w^*$ is a
local minimum of  $C$ and  that the Hessian  of $C$ at  $w^*$ is positive definite. Then
there exist $0<\delta \leq \mu $ such that 
  \begin{enumerate}
  \item $\Vert \nabla C(w)\Vert \leq \mu \Vert w-w^*\Vert $
  \item $\langle\nabla C(w), w-w^*\rangle \geq \delta \Vert w-w^*\Vert^2$
  \end{enumerate}
  for every $w$  in a neighborhood of $w^*$.
\end{remark}

\begin{proof}
Indeed, for every $\eta>0$, there exists $r>0$ such that 
$\Vert D^{2}C(w)-D^{2}C(w^{*})\Vert <\eta$ provided that   $\Vert w-w^{*}\Vert <r$ by continuity of $D^{2}C$.
Hence
  \begin{equation*}
    \begin{split}
      |\langle\nabla C(w),h\rangle|= &\left|\int_{0}^{1}
      D^{2}C(w^{*}+t(w-w^{*}))(h,w-w^{*})\,dt\right| \\ \leq &(\eta +\Vert D^{2}C(w^{*})\Vert)\,\Vert
      h\Vert \, \Vert w-w^{*}\Vert
    \end{split}
  \end{equation*}
 which implies 
  \begin{equation*}
    \Vert \nabla C(w)\Vert  \leq \ (\eta+ \Vert D^{2}C(w^{*})\Vert)\,\Vert w-w^*\Vert.
  \end{equation*}
  Thus $(1)$ holds for $\mu = (\eta+ \Vert D^{2}C(w^{*})\Vert)$.

On the other hand,
  \begin{align*}
    \langle\nabla C(w),w-w^{*}\rangle=  &\int_{0}^{1}                                          D^{2}C(w^{*}+t(w-w^{*}))(w-w^{*},w-w^{*})
                                         \,dt  \\ \geq & D^{2}C(w^{*})(w-w^{*},w-w^{*})- \eta  \Vert w-w^{*}\Vert^{2}
  \end{align*}
  provided that  $\Vert w-w^{*}\Vert <r$.
  As $D^{2}C(w^{*})$ is a bilinear form  positive definite,  
  $ N(h)=\left(D^{2}C(w^{*})(h,h)\right)^{1/2}$ is a  norm, and consequently there  exists $m>0$ such that
  \begin{equation*}
    \left(D^{2}C(w^{*})(h,h)\right)^{1/2}\geq m\Vert h\Vert
  \end{equation*}
for every  $h$. Taking $\eta=\frac12 m^{2}$ we obtain $(2)$ for $\delta =\frac{m}{2}$
 \end{proof}

 Last, we are going to give and example that provides a wide family of functions that satisfy 
 the requirements of Theorem \ref{T1}.
 
 \begin{example}\label{ejemplo}
  Let $\varphi :\mathbb{R}^+\to \mathbb{R}^+$ be a locally Lipschitz function satisfying
that  there exist
  $\delta '>0$ y $\mu '>0$ such that
  \begin{equation}\label{eq:5}
    \delta' r\leq \varphi'(r)\leq \mu' r
  \end{equation}
  for every  $r\in \mathbb{R}^+$ where $\varphi '$ exists. 
  Then for every norm
  $\Vert \ \Vert^{\star} $ 
  the function 
  \begin{equation*}
    C(w)=\varphi (\Vert w\Vert^{\star})
  \end{equation*}
 enjoy conditions 
  \ref{item:304} and \ref{item:314} of Theorem \ref{mainlocal}
 whenever $w^*=0$.
\end{example}
\begin{proof}
  It is immediate that  $C$ is locally Lipschitz.
  Also it is not difficult to see that 
  $C$ is differentiable at $w\neq 0$ if and only if
    $\varphi$ is differentiable at $\Vert w\Vert^{\star}$ and $\Vert \ \Vert^{*}$
    is differentiable at $w$. Indeed, if $C$ is differentiable at $w$, then
      \begin{equation*}
        \lim_{s\to 0}\frac{\varphi(\Vert w\Vert^{\star}+s)-\varphi(\Vert
          w\Vert^{\star})}{s}= \lim_{s\to 0}\frac{C(w+s \frac{w}{\Vert w\Vert^{\star}})-C(w)}{s}=
          \frac{1}{\Vert w\Vert^{\star}}\langle\nabla C(w),w \rangle
      \end{equation*}
      and
\begin{multline*}
  \lim_{h\to 0} \frac{\Vert w+h\Vert^{\star}-\Vert w\Vert^{\star}-\frac{1}{\varphi’(\Vert w\Vert^{\star})}\langle\nabla
            C(w),h\rangle}{\Vert h\Vert} 
  =
            \\
        \frac1{\varphi’(\Vert w\Vert^{\star})}\lim_{h\to 0}\left[\varphi’(\Vert
          w\Vert^{\star})-\frac{\varphi(\Vert
            w+h\Vert^{\star})-\varphi(\Vert w\Vert^{\star})}{\Vert
            w+h\Vert^{\star}-\Vert w\Vert^{\star}}\right]\frac{\Vert
          w+h\Vert^{\star}-\Vert w\Vert^{\star}}{\Vert h\Vert } 
          \\
          +\lim_{h\to0} \frac{\varphi(\Vert
          w+h\Vert^{\star})-\varphi(\Vert w\Vert^{\star})-\langle\nabla
          C(w),h\rangle}{\Vert h\Vert}=0.
      \end{multline*}
   Once established this equivalence, we observe first that 
    for every  $h\in\mathbb{R}^{N}$ such that $\Vert h\Vert=1$, we have
  \begin{equation*}
    \begin{split}
      \langle\nabla C(w),h\rangle=& \lim_{s\to 0^{+}}\frac{\varphi(\Vert w+sh\Vert^{\star})-\varphi(\Vert
                                    w\Vert^{\star})}{s} \\ =& \lim_{s\to 0^{+}}\frac{\varphi(\Vert w+sh\Vert^{\star})-\varphi(\Vert
                                                              w\Vert^{\star})}{\Vert w+sh\Vert^{\star}-\Vert w\Vert^{\star}}  \frac{\Vert w+sh\Vert^{\star}-\Vert                                                              w\Vert^{\star}}s\leq \varphi’(\Vert w\Vert^{\star})
    \end{split}
  \end{equation*}
  which implies
  $\Vert \nabla C(w)\Vert\leq \mu'\Vert w\Vert^{\star}$ and consequently we have that 
$\Vert \nabla C(w)\Vert \leq \mu \Vert w\Vert $ for some $\mu >0$.

For the other condition, we observe that
  \begin{equation*}
    \begin{split}
      \langle\nabla C(w),w\rangle= & \lim_{s\to 0^{+}}\frac{\varphi(\Vert w+sh\Vert^{\star})-\varphi(\Vert
      w\Vert^{\star})}{\Vert w+sw\Vert^{\star}-\Vert
      w\Vert^{\star}}  \frac{\Vert w+sw\Vert^{\star}-\Vert
      w\Vert^{\star}}{s} \\ = &  \varphi’(\Vert w\Vert^{\star})\Vert
      w\Vert^{\star} \geq  \delta’{\Vert
      w\Vert^{\star}}^{2}
    \end{split}
  \end{equation*}
  and consequently
  $\langle \nabla C(w),w\rangle \geq \delta \Vert w\Vert^{2}$ for some $\delta>0$.
\end{proof}

We finish this section with some comments on the constants that arise in the proofs.
Assume that $\delta =\mu \leq  \frac{1+2\sqrt{11}}{16}$ (this happens for instance
if $C(w)=\frac{1}{N}\Vert w\Vert_2^2$ which is an usual choice) Then we may take $A=1$. Conditions in Theorem \ref{mainlocal} reads as follows:
\begin{align*}
  \frac{\sqrt\varepsilon}{2\mu}<&
  \alpha\leq \frac{\sqrt\varepsilon}{2\mu}\left(1+\frac{\mu}{4}\right) \quad \text{and}
  \quad \beta_{1} =1-\frac{\sqrt{\varepsilon}}{2\alpha\mu},
  \\
  D=&
  \Biggl( 1- \frac{9}{16}  +
  \frac{1}{4}\Biggr) ^{\frac{1}{2}}=\frac{\sqrt{11}}{4},
\\
  L_{0}=& \max\biggl\{\beta_{1}+\mu(1-\beta_{1}),
  \beta_{2}+\mu^{2}(1-\beta_{2}),  \frac{\alpha}{\sqrt{\varepsilon}} -\frac{1}{2\mu}+
  \frac{\sqrt{11}}{4}  \biggr\} \\<&
  \max\biggl\{\beta_{1}+\mu(1-\beta_{1}),
  \beta_{2}+\mu^{2}(1-\beta_{2}), 
  \frac{1+2\sqrt{11}}{8}  \biggr\} \\ =&
  \max\biggl\{\mu +(1-\mu )\beta_{1},
  \beta_{2}+\mu^{2}(1-\beta_{2}), 
  \frac{1+2\sqrt{11}}{8}  \biggr\} 
\\ =& \max\biggl\{1-\frac{\sqrt{\varepsilon}}{2\alpha\mu}+\frac{\sqrt{\varepsilon}}{2\alpha},
  \beta_{2}+\mu^{2}(1-\beta_{2}), 
  \frac{1+2\sqrt{11}}{8}  \biggr\}.
\end{align*}
The second term is smaller than the third one provided that $\beta_2$ is small.
For the first term we have
\begin{equation*}  1-\frac{\sqrt{\varepsilon}}{2\alpha\mu}+\frac{\sqrt{\varepsilon}}{2\alpha}\leq
  1-\frac{\sqrt{\varepsilon}}{2\alpha\mu}+\mu \leq 1+\mu -(1+\mu /4)^{-1}<
  (1+\mu )(1+\mu /4)-1
\end{equation*}
which is also smaller than the third one. Hence $L_0= \frac{1+2\sqrt{11}}{8}$.

On the other hand 
\begin{align*}
  K_0=&2K_1K_2=
  2 \frac{\alpha}{\sqrt{\varepsilon}} \Big( \mu +\frac{\beta_1}{1-\beta_1}\Big)
  \Big( (1-\beta_1)+\sqrt{1-\beta_2}\Big)^{-1}
  \\
\leq&  
  2 \frac{\alpha}{\sqrt{\varepsilon}} \Big( \mu +\frac{\beta_1}{1-\beta_1}\Big) \leq
  \Big( \frac{1}{\mu}+\frac{1}{4}\Big)  \Big( \mu +\frac{\beta_1}{1-\beta_1}\Big)
  =1+\frac{\mu}{4}+
  \Big( \frac{1}{\mu}+\frac{1}{4}\Big)  \frac{\beta_1}{1-\beta_1}
\end{align*}
provided that $\beta_2$ is small enough. As
\begin{equation*}  \frac{\beta_1}{1-\beta_1}=\frac{1-\frac{\sqrt{\varepsilon}}{2\alpha\mu}}{\frac{\sqrt{\varepsilon}}{2\alpha\mu}}=
  \frac{2\alpha \mu}{\sqrt{\varepsilon}}-1\leq \frac{\mu}{4}
\end{equation*}
we obtain that 
\begin{equation*}
  K_0\leq 1+\frac{\mu}{4}+
  \Big( \frac{1}{\mu}+\frac{1}{4}\Big)  \frac{\mu}{4}=\frac{5}{4}+\frac{5}{16}\mu.
\end{equation*}
Moreover, in Theorem \ref{T2} we may also assume that $\beta =\beta_1$.

Joining all these ingredients, we deduce that in Theorem \ref{mainlocal} we may take
\begin{equation*}  L=\frac{L_0+1}{2}=\frac{\cfrac{1+2\sqrt{11}}{8}+1}{2}=\frac{9+2\sqrt{11}}{16}
\end{equation*}
and $n_0=1$ since 
\begin{equation*}
  K_0\beta \leq \frac{5}{4}\beta_1+\frac{5}{16}\mu \beta_1\leq \frac{25}{16}\beta_1\leq
  \frac{1-  \cfrac{1+2\sqrt{11}}{8}}{2}
  =\frac{1-L_0}{2}=L-L_0
\end{equation*}
if $\beta_1\leq \frac{8}{25}(1-  \frac{1+2\sqrt{11}}{8})=\frac{7-2\sqrt{11}}{25}$
for instance. Therefore we may take
\begin{equation*}
  K\leq K_0+L_0\leq \frac{5}{4}+\frac{5}{16}\mu+  \frac{1+2\sqrt{11}}{8}<
  \frac{27+4\sqrt{11}}{16}=1+\frac{11+4\sqrt{11}}{16}<3.
\end{equation*}

\section{Approximation to Convergence Basin}

 Theorem $3.4$  in \cite{DMU} provides conditions that guarantee that the algorithm, in a different form 
that we show,  reaches
the convergence basin. The aim of this section is to prove that a similar result holds also in our setting.

Assume that $C$ is  a Lipschitz function.
We will define recursively $m_n$, $v_n$ and $w_n$ as in generalized Adam.
We start with a serie of Lemmas that provide us  estimations that we will require throughout this section.
Observe that these estimations are independent of the parameters choice.
 We denote $\displaystyle \sigma_{n}=\max_{i=1,\dots ,n}\Vert \zeta_{w_i}\Vert$.

\begin{lemma}\label{sec:aprox-la-cuenca-1}
  If $m_0=0$ then
  \begin{equation*}
    \biggl\Vert \frac{m_{n+1}}{\sqrt{v_{n+1}+(\varepsilon ,\dots ,\varepsilon )}}\biggr\Vert \leq
    \frac{(1-\beta_1^{n+1})}{\sqrt{\varepsilon}}\sigma_{n}.
  \end{equation*}
\end{lemma}
\begin{proof}
  From
  \begin{equation*}
    m_{n+1}=\beta_1 ^{n+1}m_0+(1-\beta_1 )\sum_{i=0}^n \beta_1 ^{n-i}g_i=
    (1-\beta_1 )\sum_{i=0}^n \beta_1 ^{n-i}\zeta_{w_i},
  \end{equation*}
  (see Remark \ref{formulasbasicas}), we deduce
  \begin{equation*}
    \Vert m_{n+1}\Vert \leq
    (1-\beta_1 )\sum_{i=0}^n\beta_1 ^{n-i}\sigma_{n}=
    (1-\beta_1^{n+1})\sigma_{n}.
  \end{equation*}
\end{proof}

\begin{lemma}\label{sec:aprox-la-cuenca-2}
 If  $m_0=v_0=0$, then
  \begin{equation*}
    \biggl\langle
    \zeta_{w_n},\frac{m_{n+1}}{\sqrt{v_{n+1}+(\varepsilon ,\dots ,\varepsilon )}}
    \biggr\rangle \geq
    \Vert \zeta_{w_n}\Vert^2\bigg(
    \frac{(1-\beta_1)}
    {\sigma_n \sqrt{(1-\beta_2^{n+1})}+\sqrt{\varepsilon}}-
    \frac{(\beta_1-\beta_1^{n+1})\sigma_n}{\Vert \zeta_{w_n}\Vert \sqrt{\varepsilon}}
    \bigg)
  \end{equation*}
  for every $n$.
\end{lemma}
\begin{proof}
  \allowdisplaybreaks
  \begin{align*}
      &\biggl\langle \zeta_{w_n},\frac{m_{n+1}}{\sqrt{v_{n+1}+(\varepsilon ,\dots ,\varepsilon )}}
      \biggr\rangle =  \biggl\langle \frac{\zeta_{w_n}}{\sqrt{v_{n+1}+(\varepsilon ,\dots ,\varepsilon )}},m_{n+1}
      \biggr\rangle  \\ =&
      \biggl\langle \frac{\zeta_{w_n}}{\sqrt{v_{n+1}+(\varepsilon ,\dots ,\varepsilon )}},\beta_1m_{n}+(1-\beta_1)\zeta_{w_n}
      \biggr\rangle \\= &
      \beta_1  \biggl\langle \frac{\zeta_{w_n}}{\sqrt{v_{n+1}+(\varepsilon ,\dots ,\varepsilon )}},m_{n}
      \biggr\rangle +(1-\beta_1) \biggl\langle \frac{\zeta_{w_n}}{\sqrt{v_{n+1}+(\varepsilon ,\dots ,\varepsilon )}},\zeta_{w_n}
      \biggr\rangle \\ \geq&
      (1-\beta_1) \Vert \zeta_{w_n}\Vert^2 \frac{1}{\sqrt{(1-\beta_2^{n+1})\sigma_n^2+\varepsilon}}+
      \beta_1  \biggl\langle \frac{\zeta_{w_n}}{\sqrt{v_{n+1}+(\varepsilon ,\dots ,\varepsilon )}},m_{n}
      \biggr\rangle \\ \geq&
      (1-\beta_1) \Vert \zeta_{w_n}\Vert^2 \frac{1}{\sqrt{(1-\beta_2^{n+1})\sigma_n^2+\varepsilon}}-
      \frac{\beta_1-\beta_1^{n+1}}{\sqrt{\varepsilon}}
                             \Vert \zeta_{w_n}\Vert \sigma_n \\ \geq &
      \Vert \zeta_{w_n}\Vert^2\bigg(
      \frac{(1-\beta_1)}
      {\sigma_n \sqrt{(1-\beta_2^{n+1})}+\sqrt{\varepsilon}}-
      \frac{(\beta_1-\beta_1^{n+1})\sigma_n}{\Vert \zeta_{w_n}\Vert \sqrt{\varepsilon}}
      \bigg)
  \end{align*}
by Remark \ref{formulasbasicas} and Lemma \ref{sec:aprox-la-cuenca-1}.
\end{proof}

Observe that if $C$ is $\sigma$-Lipschitz, then $\Vert \zeta_{w_n}\Vert \leq \sigma$ for every
$n$, hence we obtain the following inequality:

\begin{remark}\label{sec:aprox-la-cuenca-3}
If If  $m_0=v_0=0$
and $\Vert \zeta_{w_n}\Vert \leq \sigma$ for every
$n$
 \begin{equation*}
    \biggl\langle
    \zeta_{w_n},\frac{m_{n+1}}{\sqrt{v_{n+1}+(\varepsilon ,\dots ,\varepsilon )}}
    \biggr\rangle \geq
    \Vert \zeta_{w_n}\Vert^2\bigg(
    \frac{(1-\beta_1)}
    {\sigma \sqrt{(1-\beta_2^{n+1})}+\sqrt{\varepsilon}}-
    \frac{(\beta_1-\beta_1^{n+1})\sigma}{\Vert \zeta_{w_n}\Vert \sqrt{\varepsilon}}
    \bigg).
  \end{equation*}
\end{remark}

\begin{lemma}\label{sec:aprox-la-cuenca-5}
  Let  us set $\sigma >0$. 
  If $\eta \in (0,\sigma )$, $\beta_1 \in (0,\frac{\eta}{\eta +\sigma})$, 
   $\varepsilon$ satisfies
   \begin{equation*}
     \sqrt{\varepsilon}\Big( \frac{(1-\beta_1)\eta}{\sigma \beta_1}-1\Big) >\sigma,
   \end{equation*}
  and $0<\beta_{2}<1$, then there exist 
  $\theta_1,\theta_2>0$ such that
  \begin{equation*}
    \bigg(
    \frac{(1-\beta_1)}
    {\sigma \sqrt{(1-\beta_2^{n+1})}+\sqrt{\varepsilon}}-
    \frac{(\beta_1-\beta_1^{n+1})\sigma}{\eta \sqrt{\varepsilon}}
    \bigg) \geq
    \frac{(\beta_1-\beta_1^{n+1})
      \theta_1 \theta_2}
    {\sqrt{\varepsilon}(\sigma +\sqrt{\varepsilon})}
  \end{equation*}
for every $n$.
\end{lemma}
\begin{proof}
 We define
  \begin{equation*}
    \theta_1:=\frac{(1-\beta_1)\eta}{\sigma \beta_1}-1>
    \frac{(1-\frac{\eta}{\eta +\sigma})\eta}{\sigma \frac{\eta}{\eta +\sigma}}-1=0
    \quad\text{and }\quad
    \theta_2:=\sqrt{\varepsilon}-\frac{\sigma}{\theta_1}>0.
  \end{equation*}
Then
  \begin{align*}
    &\frac{(1-\beta_1)}
      {\sigma \sqrt{(1-\beta_2^{n+1})}+\sqrt{\varepsilon}} -
      \frac{(\beta_1-\beta_1^{n+1})\sigma}{ \eta
      \sqrt{\varepsilon}} \\=
    &
      \frac{(1-\beta_1)\eta \sqrt{\varepsilon}-
      (\sigma \sqrt{(1-\beta_2^{n+1})}+\sqrt{\varepsilon})(\beta_1-\beta_1^{n+1})\sigma}
      { \eta \sqrt{\varepsilon}(\sigma
      \sqrt{(1-\beta_2^{n+1})}+\sqrt{\varepsilon})}
    \\=
    &
      \sigma (\beta_1-\beta_1^{n+1})
      \frac{\sqrt{\varepsilon}
      \Big( \frac{(1-\beta_1) \eta }{\sigma (\beta_1-\beta_1^{n+1})}-1\Big) -
      \sigma \sqrt{(1-\beta_2^{n+1})}}
      {\eta \sqrt{\varepsilon}(\sigma    \sqrt{(1-\beta_2^{n+1})}+\sqrt{\varepsilon})}
    \\ \geq
    &
      \sigma (\beta_1-\beta_1^{n+1})
      \frac{\sqrt{\varepsilon}\Big( \frac{(1-\beta_1)\eta}{\sigma \beta_1}-1\Big)
      -\sigma \sqrt{(1-\beta_2^{n+1})}}
      {\eta \sqrt{\varepsilon}(\sigma
      \sqrt{(1-\beta_2^{n+1})}+\sqrt{\varepsilon})}
    \\ =
    &
      \sigma (\beta_1-\beta_1^{n+1})\theta_1
      \frac{\sqrt{\varepsilon} -
      \frac{\sigma \sqrt{(1-\beta_2^{n+1})}}
      {\theta_1}}
      {\eta \sqrt{\varepsilon}(\sigma   \sqrt{(1-\beta_2^{n+1})}+\sqrt{\varepsilon})}
    \\ \geq
    &
      \sigma (\beta_1-\beta_1^{n+1})
      \theta_1
      \frac{\sqrt{\varepsilon} -
      \frac{\sigma}{\theta_1}}
      { \eta \sqrt{\varepsilon}(\sigma
      +\sqrt{\varepsilon})} \\ =
    &
      \sigma (\beta_1-\beta_1^{n+1})
      \theta_1
      \frac{\theta_2}
      {\eta \sqrt{\varepsilon}(\sigma +\sqrt{\varepsilon})}\geq
      \frac{(\beta_1-\beta_1^{n+1})
      \theta_1 \theta_2}
      {\sqrt{\varepsilon}(\sigma +\sqrt{\varepsilon})}.
  \end{align*}
\end{proof}

In order to achieve the main result of this section, we must require an additional
condition on the function $C$.  Namely:
there exist  $M,R>0$ such that
    \begin{equation}
      C(w')\leq C(w)+\langle \nabla C(w),w'-w\rangle +\frac{M}{2}\Vert
      w'-w\Vert^2\label{eq:3}
      \tag{$\star$}
    \end{equation}
     provided that  $\Vert w-w'\Vert <R_0$ and that $C$ is differentiable at  $w$.

 Observe that \eqref{eq:3} implies that
  \begin{equation}
    C(w')\leq C(w)+\langle \zeta_w,w'-w\rangle +\frac{M}{2}\Vert w'-w\Vert
    ^2\label{eq:4}
    \tag{$\star\star$}
  \end{equation}
 provided that  $\Vert w-w'\Vert <R_0$. 
We will prove that this condition allow us to control the behavior of function $C$ on the sequence $(w_n)_n$.

First we precise the parameter in Generalized Adam as follows: 
\begin{equation*}
    \alpha_n=
    \frac{\Bigl\langle \zeta_{w_n},\frac{m_{n+1}}{\sqrt{v_{n+1}+(\varepsilon ,\dots ,\varepsilon)}}\Bigr\rangle}
    {M \Bigl\Vert \frac{m_{n+1}}{\sqrt{v_{n+1}+(\varepsilon ,\dots ,\varepsilon )}}\Bigr\Vert^2}.
  \end{equation*}
Moreover we will start the algorithm with $(m_{0},v_{0},w_{0})=(0,0,w_{0})$. Observe that
if $C$ is $\sigma$-Lipschitz and we take the parameters as in Lemma \ref{sec:aprox-la-cuenca-5}
then $\alpha_n>0$.

From these preliminar estimations we may deduce the main theorem of this section.

\begin{theorem}\label{maincuenca}
  Assume that $C:\mathbb{R}^n\to \mathbb{R}^+$ is a $\sigma$-Lipschitz function that satisfies
  ~\eqref{eq:3}. Assume also that $w^*$ is a global minimum of $C$. 
  
 Then for every
$\eta >0$, $\beta_1$, $\beta_2$, and $\varepsilon$
  satisfying the restrictions  of Lemma \ref{sec:aprox-la-cuenca-5}
there exists $(\alpha_n)_n$ such that
  $w_n=w^*$ eventually  or there exists $n_0$ such that
  $0<\Vert \zeta_{w_{n_0}}\Vert \leq \eta$.
\end{theorem}
\begin{proof}
 We start the generalized Adam with $(m_{0},v_{0},w_{0})=(0,0,w_{0})$,
  defining  $\alpha_n$ as above.
 Assume that there exist  $\eta >0$,
  $\beta_1$, $\beta_2$ and $\varepsilon$ satisfying the restrictions of Lemma \ref{sec:aprox-la-cuenca-5}
  such that 
  $\Vert \zeta_{w_{n}}\Vert >\eta$ for every  $n$.  
  First, we observe that
  \begin{equation*}
    \Vert w_{n+1}-w_n\Vert \leq \frac{\Bigl\langle \zeta_{w_n},\frac{m_{n+1}}{\sqrt{v_{n+1}+(\varepsilon ,\dots ,\varepsilon)}}\Bigr\rangle}
    {M \Bigl\Vert \frac{m_{n+1}}{\sqrt{v_{n+1}+(\varepsilon ,\dots ,\varepsilon )}}\Bigr\Vert} \leq \frac{\sigma}{M}<R_0
  \end{equation*}
  increasing $M$ if necessary. 
  
  Then invoking 
   \eqref{eq:4}, we obtain
  \begin{equation*}
    C(w_{n+1})\leq C(w_n)+\langle \zeta_{w_n},w_{n+1}-w_n\rangle +\frac{M}{2}\Vert w_{n+1}-w_n\Vert^2.
  \end{equation*}
Hence
 \begin{equation*}
    \begin{split}
      C(w_{n+1})&-C(w_n) \\ &\leq
                          \alpha_n \bigg(
                          -\Big \langle \zeta_{w_n},\frac{m_{n+1}}{\sqrt{v_{n+1}+(\varepsilon ,\dots ,\varepsilon )}}\Big \rangle  +\frac{M\alpha_n}{2}\Big| \Big| \frac{m_{n+1}}{\sqrt{v_{n+1}+(\varepsilon ,\dots ,\varepsilon )}}\Big| \Big| ^2
                          \bigg) \\ & \leq
                                      -\frac{1}{2M} \frac{\Big \langle \zeta_{w_n},\frac{m_{n+1}}{\sqrt{v_{n+1}+(\varepsilon ,\dots ,\varepsilon )}}\Big \rangle ^2}
                                      {\Bigl\Vert \frac{m_{n+1}}{\sqrt{v_{n+1}+(\varepsilon ,\dots ,\varepsilon )}}\Bigr\Vert^2}.
    \end{split}
  \end{equation*}

 Therefore, by Lemma \ref{sec:aprox-la-cuenca-1}, Remark \ref{sec:aprox-la-cuenca-3} and
 Lemma \ref{sec:aprox-la-cuenca-5}, we deduce
  \begin{equation*}
    \begin{split}
      C(w_{n+1})-C(w_n)\leq &
                              -\frac{\Vert \zeta_{w_n}\Vert ^4\bigg(
                              \frac{(1-\beta_1)}
                              {\sigma \sqrt{(1-\beta_2^{n+1})}+\sqrt{\varepsilon}}-
                              \frac{(\beta_1-\beta_1^{n+1})\sigma}{\Vert \zeta_{w_n}\Vert \sqrt{\varepsilon}}
                              \bigg) ^2
                              }{2M\frac{(1-\beta_1^{n+1})^2\sigma ^2}{\varepsilon}} \\ \leq
                            &
                              -\frac{\Vert \zeta_{w_n}\Vert ^4\bigg(
                              \frac{\beta_1(1-\beta_1^n)
                              \theta_1 \theta_2}
                              {\sqrt{\varepsilon}(\sigma +\sqrt{\varepsilon})}
                              \bigg) ^2
                              }{2M \frac{(1-\beta_1^{n+1})^2\sigma ^2}{\varepsilon}}=
                              -\frac{\Vert \zeta_{w_n}\Vert ^4\bigg(
                              \frac{\beta_1(1-\beta_1^n)
                              \theta_1 \theta_2}
                              {(\sigma +\sqrt{\varepsilon})}
                              \bigg) ^2
                              }{2M (1-\beta_1^{n+1})^2\sigma ^2} \\\leq
                            &
                              -\frac{\eta^4\bigg(
                              \frac{\beta_1(1-\beta_1)
                              \theta_1 \theta_2}
                              {(\sigma +\sqrt{\varepsilon})}
                              \bigg) ^2
                              }{2M \sigma ^2}.
    \end{split}
  \end{equation*}
  This implies that $w_n=w^*$ eventually. Alternatively, there exists $n$ such that $\Vert \zeta_{w_n}\Vert \leq \eta$.
\end{proof}

\begin{remark}\label{remT}
Observe that we may estimate the number of steps that we require in order to attain the minimum
or alternatively to generate a generalized gradient $\zeta_w$  such that 
$\Vert \zeta_w\Vert \leq \eta$.
 Indeed, if we denote 
  \begin{equation*}
    \begin{split}
      s=&\frac{\eta^4\bigg(
          \frac{\beta_1(1-\beta_1)
          \theta_1 \theta_2}
          {(\sigma +\sqrt{\varepsilon})}
          \bigg) ^2
          }{2M\sigma ^2}=
          \frac{\eta^4\bigg(
          \frac{(1-\beta_1)
          ((1-\beta_1)\eta -\sigma \beta_1) \sqrt{\varepsilon}}
          {\sigma (\sigma +\sqrt{\varepsilon})}
          -
          \frac{\beta_1(1-\beta_1)
          \sigma}
          {(\sigma +\sqrt{\varepsilon})}
          \bigg) ^2
          }{2M\sigma ^2} \\ = &
                                \frac{\eta^4\bigg(
                                (1-\beta_1)
                                ((1-\beta_1)\eta -\sigma \beta_1) \sqrt{\varepsilon}
                                -
                                \beta_1(1-\beta_1)\sigma ^2
                                \bigg) ^2
                                }{2M\sigma ^4(\sigma +\sqrt{\varepsilon})^2} \\= &
                                                                                   \frac{\varepsilon \eta^4\bigg(
                                                                                   (1-\beta_1)
                                                                                   ((1-\beta_1)\eta -(\sigma +\frac{\sigma^2}{\sqrt{\varepsilon}}) \beta_1)
                                                                                   \bigg) ^2
                                                                                   }{2M\sigma ^4(\sigma +\sqrt{\varepsilon})^2}
    \end{split}
  \end{equation*}
 and $n_0=\big[ \frac{C(w_0)}{s}\big]$, then $\Vert \zeta_{w_n}\Vert \leq \eta$ for some $n\leq n_0$.
\end{remark}

Looking back at the Example \ref{ejemplo}, for a function defined as
$C(w)=\varphi (\Vert w\Vert^{\star})$, we have the following example.

\begin{example}\label{ejemplobis}
Let $R_0,M_0,M_1$ be positive constants
  Let $\varphi :\mathbb{R}^+\to \mathbb{R}^+$ be a  Lipschitz function satisfying
  $\varphi '(r)\geq 0$ and 
  \begin{equation*}
    \varphi (r_1)-\varphi (r)-\varphi '(r)(r_1-r)\leq M_0(r_1-r)^2
  \end{equation*}
 provided that $|r-r'|<R_0$, whenever $\varphi'(r)$ exists. 
  Assume also that the norm $\Vert \ \Vert^*$ satisfies
   \begin{equation}
      \label{eq:1} \tag{$\star\star\star$}
      \Vert w’\Vert^{\star}-\Vert
      w\Vert^{\star}-\langle\nabla\Vert\ \Vert^{\star}(w),w’-w\rangle \leq
      M_{1}\Vert w’-w\Vert^{2}
    \end{equation}
 for every $w$ where $\Vert \ \Vert^*$ is differentiable.  
Then the function 
  \begin{equation*}
    C(w)=\varphi (\Vert w\Vert^{\star})
  \end{equation*}
 enjoys  condition \eqref{eq:3}.
\end{example}
\begin{proof}
Let us see it.
 \begin{align*}
      C(w_1)& -C(w)-\langle \nabla C(w),w_1-w\rangle \\ =&
                    \varphi (\Vert w_1\Vert ^{\star})-\varphi (\Vert w\Vert
                    ^{\star})-\varphi'(\Vert w\Vert
                    ^{\star})\langle\nabla\Vert \ \Vert
                    ^{\star}(w),w_1-w\rangle 
      \\ \leq
                  &
                    \varphi (\Vert w_1\Vert^{\star})-\varphi (\Vert w\Vert^{\star})- \varphi'(\Vert w
                    \Vert^{\star})(\Vert w_1\Vert^{\star}-\Vert
                    w\Vert^{\star}) + M_{1}\varphi'(\Vert w
                    \Vert^{\star}) \Vert w_1-w\Vert^2\\ \leq
                  &
                    M_0(\Vert w_1\Vert^{\star}-\Vert w\Vert^{\star})^2 + M_{1}\varphi'(\Vert w
                    \Vert^{\star}) \Vert w_1-w\Vert^2       \\    \leq &
                                                                         M_0(\Vert
                                                                         w_1-w\Vert^{\star})^2 + M_{1}\varphi'(\Vert w
                                                                         \Vert^{\star}) \Vert w_1-w\Vert^2 
                                                                         \leq  \frac{M}{2}\Vert w_1-w\Vert^2
    \end{align*}   
provided that $\Vert w_1-w\Vert <\frac{R_0}{A_0}$, since 
$\big| \Vert w_1\Vert^*- \Vert w\Vert^* \big| \leq A_0\Vert w_1-w\Vert $ for a suitable constant $A_0$.
\end{proof}

Observe that condition \eqref{eq:1} is satisfied if
$\Vert \ \Vert^*$ is locally linear whenever it is differentiable, since in this situation
$D(\Vert \ \Vert^{\star})(w)=\Vert \ \Vert^{\star}$, and consequently
\begin{equation*}
  \Vert w'\Vert^{\star}-\Vert w \Vert^{\star}-\Vert w'-w \Vert^{\star}\leq 0.
\end{equation*}
Incidentally, this condition happens if and only if the unit ball
is the convex hull of a finite number of points. 

Another interesting situation where (\ref{eq:1}) is satisfied  happens when
$\Vert \ \Vert^{\star}$ is an euclidean norm since in this case
\begin{equation*}
  D^2(\Vert \ \Vert^{\star})(w)=\frac{1}{\Vert w\Vert ^{\star}}I-\frac{1}{(\Vert w \Vert^{\star} )^3}
  \langle w, \cdot \rangle ^{\star}
\end{equation*}
which is bounded away from $0$ (we only require that).

\section{Global convergence}

Putting together the main theorems of the two previous sections we obtain the following
global result.

\begin{theorem}\label{mainglobal}
 Assume that $C:\mathbb{R}^n\to \mathbb{R}^+$ is  $\sigma$-Lipschitz,
 that satisfies ~\eqref{eq:3}, and that attains a unique global minimum at 
 $w^*\in \mathbb{R^N}$. Assume also that there exist
 $R\in (0,1)$, and $\delta \leq \mu$, such that whenever  $C$ is differentiable at $w$
  the following conditions hold:
  \begin{enumerate}[label=\textup{(\arabic*)}]
  \item \label{item:18}  $\Vert \nabla C(w)\Vert \leq \mu \Vert w-w^*\Vert $ for every $w\in B(w^*,R)$,
  \item \label{item:19} 
    $\displaystyle\langle \nabla C(w),w-w^*\rangle \geq
    \delta \Vert w-w^*\Vert^2$ for every $w\in B(w^*,R)$.
  \end{enumerate}
If  $w_0\in B(w^*,R)$ then there exist suitable parameters 
  $\beta_1^*,\beta_2^*,\varepsilon^*, (\alpha_n)_n$ and a step $n_0$
 such that at some step $n\leq n_0$, Generalized Adam either attains $w^*$
 or allows us to define a new algorithm with parameters
  $\beta_1,\beta_2,\varepsilon , \alpha$ such that  $(w_n)_n$ converges
  exponentially to $w^*$.
\end{theorem}

\begin{proof}
We choose parameters 
 $A$, $\varepsilon$, $\alpha$, and $\beta_1\beta_2\in (0,1)$ satisfying
   \begin{enumerate}
  \item \label{item:20}
$A>\mu$ (If $\mu <1$ then we choose  $A=1$),
    \item \label{item:21}
    $\displaystyle 
    \frac{\delta\sqrt\varepsilon}{2\mu^{2}}<
      \alpha\leq \frac{\delta\sqrt\varepsilon}{2\mu^{2}}\left(1+\frac{\delta}{4A}\right) 
      $ 
    \end{enumerate}
    and
    \begin{enumerate}[resume*]
  \item\label{item:22}
    $\displaystyle
     \beta_1 =1-\frac{\delta \sqrt{\varepsilon}}{2\alpha\mu^2}$.
  \end{enumerate}

Then we put $r=\min (1,R,\mu^{-1})$ and $\eta =\delta \frac{r}{K}$, where
  $K:=K(L_0,\beta_1,\beta_2,\varepsilon ,\alpha )$ is defined as in  Theorem \ref{mainlocal}.

  We start Generalized Adam  with $m_0=v_0=0$,
  $w_0\in \mathbb{R}^N$, and parameters
  $0<\beta_1^*<\frac{\eta}{\eta +\sigma}$, $0<\beta_2^*<1$, $\varepsilon^*$
  such that
  \begin{equation*}
    \sqrt{\varepsilon^*}\Big( \frac{(1-\beta_1^*)\eta}{\sigma \beta_1^*}-1\Big) >\sigma
  \end{equation*}
  and $\alpha_n$,
  as in Theorem \ref{maincuenca}. From that Theorem we deduce that 
  there exists
  $n\leq n_0$, where $n_0:=n_0(M,\eta ,\sigma , \beta_1,\varepsilon )$ is defined as in  Remark \ref{remT},
  such that $w_{n}=w^*$ or $\Vert \zeta_{w_n}\Vert \leq \eta$.

Then we start Adam with $m_{0}=v_{0}=0$ y $w_0=w_{n}$ and the parameters
$\varepsilon$, $\beta_1$,
  $\beta_2$ and $\alpha$. As
   \begin{equation*}
    \delta \Vert x_0-x^*\Vert =\delta \Vert w_0-w^*\Vert \leq
    \Vert \zeta_w\Vert   \leq  \eta =\delta \frac{r}{K}
  \end{equation*}
  we deduce $\Vert x_0-x^*\Vert \leq \frac{r}{K}$.
  Finally, we obtain the result  invoking  Theorem \ref{mainlocal}.
\end{proof}


\begin{thebibliography}{8}

\bibitem{BW} S. Bock, M. Wei{\ss}, A Proof of Local Convergence for
  the Adam Optimizer in:
  \emph{2019 International Joint Conference on Neural Networks (IJCNN)}, (2019),  1-8.
  
  \bibitem{C} F. H. Clarke, Yu. S. Ledyaev, R. J. Stern,
    R. R. Wolenski, \emph{Nonsmooth Analysis and Control Theory},
    Springer, New York, 1997.
  
  \bibitem{DMU} S. De, A. Mukherjee, E. Ullah,
  Convergence guarantees for RMSProp and ADAM in non-convex opti-
  mization and their comparison to Nesterov acceleration, preprint 
  (2018), https://doi.org/10.48550/arXiv.1807.06766.
  
\bibitem{DBB} A. Défossez, L. Bottou, F.R.  Bach, N. Usunier,  A Simple Convergence Proof of Adam and Adagrad, \emph{Trans. Mach. Learn. Res.}, (2022).
  
\bibitem{F}
  J. Ferrera,\emph{ An introduction to nonsmooth analysis}, Elsevier/Academic Press,
  Amsterdam, 2014, http://dx.doi.org/10.1016/B978-0-12-800731-0.00001-1.

\bibitem{KB} D. Kingma, J. Ba,   Adam: A Method for Stochastic
  Optimization, in:  
\emph{In International Conference on Learning Representations} (2015).

\bibitem{MR}  D. Mart\'{\i}nez Rubio, \emph{Convergence Analysis of an
  Adaptive Method of Gradient Descent}, Ph.D. thesis for the degree of MSc on Mathematics and
  Foundations of Computer Science. Wadham College, University of Oxford, 2017.
  
\bibitem{RKK}  S. J. Reddi, S. Kale,  S. Kumar, 
On the convergence of Adam and beyond, in: \emph{International Conference on Learning Representations}, 2018.
 
\end{thebibliography}
\end{document}